\documentclass{fic-l}
\usepackage{amssymb}
\newtheorem{theorem}{Theorem}[section]
\newtheorem{corollary}[theorem]{Corollary}
\newtheorem{lemma}[theorem]{Lemma}
\newtheorem{proposition}[theorem]{Proposition}

\theoremstyle{definition}
\newtheorem{definition}[theorem]{Definition}

\theoremstyle{remark}
\newtheorem{remark}[theorem]{Remark}

\numberwithin{equation}{section}

\newcommand{\PP}{{\mathbb P }}
\newcommand{\QQ}{{\mathbb Q }}

\newcommand{\ZZ}{{\mathbb Z }}
\newcommand{\oA}{{\overline{A}}}
\newcommand{\ok}{{\overline{k}}}

\newcommand{\Ku}{{K^\mathrm{u}}}
\newcommand{\oW}{{\overline{W}}}
\newcommand{\oH}{\overline{H}}
\newcommand{\oU}{\overline{U}}
\newcommand{\oV}{\overline{V}}
\newcommand{\tf}{\tilde{f}}
\newcommand{\OM}{\Omega}
\newcommand{\ep}{\varepsilon}
\newcommand{\vf}{\varphi}
\newcommand{\Fil}{\mathrm{Fil}}
\newcommand{\gr}{\mathrm{gr}}
\newcommand{\rank}{\mathrm{rank}}
\newcommand{\cT}{\mathcal{T}}
\newcommand{\cS}{\mathcal{S}}

\begin{document}

\title{Ordinary Calabi-Yau-3 Crystals }

\author{Jan Stienstra}
\address{Mathematisch Instituut, Universiteit Utrecht\\
Postbus 80.010,  3508 TA Utrecht, the Netherlands}
\email{stien@math.uu.nl}
\begin{abstract}
We show that crystals with the properties of crystalline cohomology of ordinary
Calabi-Yau threefolds in characteristic $p>0$, exhibit a remarkable similarity
with the well known structure on the cohomology of complex Calabi-Yau
threefolds near a boundary
point of the moduli space with maximal unipotent local monodromy. In
particular, there are canonical coordinates and an analogue of the prepotential
of the Yukawa coupling.
Moreover in Formulas (\ref{eq:qprops}) and (\ref{eq:yukinteger})
we show $p$-adic analogues of the integrality properties for the canonical
coordinates and the prepotential
of the Yukawa coupling, which have been observed in the examples of Mirror
Symmetry.
\end{abstract}

\maketitle

\section*{Introduction}
Calabi-Yau manifolds of dimensions $1$ and $2$, i.e. elliptic curves and
K3-surfaces, have a long and successful tradition in geometry and number
theory.
In the 1980's, in connection with developments in string theory, Calabi-Yau
manifolds of dimension $3$ moved to the forefront. Emphasis has been on
their complex and symplectic geometry, and in
particular on the variation of the Hodge structure for Calabi-Yau
threefolds near the so-called large complex structure limit.
This is the first of a series of papers in which we want to describe certain
aspects of the
arithmetic geometry of families of
ordinary Calabi-Yau threefolds
and analogies with complex Calabi-Yau threefolds near the
large complex structure limit. More specifically,
this paper discusses the associated crystals, i.e.
modules with an integrable connection, like the Gauss-Manin connection on the
cohomology of a family of Calabi-Yau threefolds.
Crystals which originate from geometry come with a Hodge filtration,
which is not preserved by the connection but instead satisfies Griffiths
transversality. Crystals which originate from geometry
in characteristic $p$ additionally carry an action of Frobenius operators
behaving in an appropriate way with respect to the connection and the Hodge
filtration: they are divisible Hodge $F$-crystals. For so-called divisible
\emph{ordinary} Hodge
$F$-crystals the space of flat sections decomposes into eigenspaces
for the Frobenius operators. The position of the Hodge filtration with respect
to this eigenspace decomposition is the
source for the canonical coordinates and the Yukawa coupling:
see Theorems \ref{thmstructure}, \ref{cy3structure}, \ref{integralities}.
The ordinariness requirement implies, exactly as the requirement of maximally
unipotent local monodromy
in the traditional analysis of complex Calabi-Yau
threefolds near the large complex structure limit \cite{M}, that there is a
\emph{filtration which is opposite to the Hodge filtration and
is invariant under the Gauss-Manin connection and whose associated graded
is of Hodge-Tate type}
(cf. \cite{D1} prop.1.3.6, \cite{D2} \S 6).

The theory of ordinary Hodge $F$-crystals has its origins in Dwork's work on
the variation of the zeta function in a family of varieties over a finite field
\cite{K1} and in the work of Serre and Tate on the formal moduli of
deformations of an ordinary elliptic curve over a field of positive
characteristic. In \cite {D1} Deligne and Illusie present
a general
theory of ordinary Hodge $F$-crystals and apply it to investigate the formal
moduli of ordinary abelian varieties and ordinary K3-surfaces.
The general theory works equally well for ordinary CY3
crystals, i.e. crystals of the kind
that arises as the crystalline cohomology of ordinary Calabi-Yau threefolds in
characteristic $p>0$. Thus, exactly as in \cite{D1}, there are canonical
coordinates on the formal moduli space. In Section \ref{section 1} we briefly
recall the general theory of \cite{D1}.
In Section 2 we take the same algebraic path as Bryant,
Griffiths, Morrison and Deligne \cite{BG,F,M,D2} in their analysis of the
variation of the Hodge structure of complex Calabi-Yau threefolds
near a maximally degenerate ($\approx$ maximally unipotent local
monodromy $\approx$ large complex structure) boundary point of the moduli
space. In particular
we find for ordinary CY3 crystals an exact analog of the cubic
form of Bryant-Griffiths ($\approx$ Yukawa coupling).
Interestingly, besides the standard characterization (\ref{eq:yukawa})
of the prepotential of this cubic form we find a characterization
(\ref{eq:yukinnerproduct})
by means of the action of the canonical lift
of Frobenius. The latter characterization remains valid if there are no
parameters, i.e. if $n=0$ in (\ref{eq:A}), while in that case
the former description is vacuous. Also for the canonical coordinates the
formalisms in the ordinary and the large complex structure situations match
very well, but the match is not perfect: in the ordinary case the canonical
coordinates take the value $1$ at the origin of the deformation space (cf.
(\ref{eq:canorig})), whereas the canonical coordinates at the large complex
structure limit point vanish.

The standard algorithms (cf. \cite{C,M}) for computing canonical coordinates
and Yukawa coupling for complex Calabi-Yau threefolds
near the ``maximally unipotent local monodromy boundary point'' use special
solutions of the Picard-Fuchs equations
associated with the family of CY3's and a nowhere vanishing global $3$-form. In
Remark \ref{Picard-Fuchs} we point out that in the ordinary case one has
basically the same relation between canonical coordinates, Yukawa coupling and
special solutions of the Picard-Fuchs equations
(although explicit computations seem here out of reach).

The present paper is a survey of structures none which -- except maybe
(\ref{eq:yukinnerproduct}) and (\ref{eq:yukinteger}) -- is new by itself,
but which appear in the literature in different contexts. So we hope
that there is an inspiring and cross-fertilizing effect from putting
them into one (con)text. For instance, being alerted that for the deformation
theory of ordinary Calabi-Yau threefolds the cubic form of Bryant-Griffiths
(Yukawa coupling) should be relevant, one may wonder how that shows up in
connection with the Hodge-Tate decomposition of $p$-adic \'etale cohomology as
in \cite{BK} p. 109.

The original motivation for our work was to find a general method of proof for
the integrality conjectures in Mirror Symmetry by reducing them to known
results or easy to check conditions on the crystalline cohomology for families
of Calabi-Yau threefolds in positive characteristics:
It has been observed in many examples in the literature on Mirror Symmetry for
Calabi-Yau threefolds, that the coefficients in
appropriate expansions of the canonical coordinates and of
the prepotential of the Yukawa coupling are integers. In the famous example of
the quintics in $\PP^4$ (see \cite{C,M}) one has for instance
$$
q =  t + 770\,t^2 + 1014275\,t^3 + 1703916750\,t^4
+ 3286569025625\,t^5+ \ldots
$$
where $t=\psi^{-5}$ is the parameter and $q$ is the canonical coordinate of the
family
$$
X_0^5+X_1^5+X_2^5+X_3^5+X_4^5-\psi X_0X_1X_2X_3X_4=0\,;
$$
the Yukawa coupling $Y$ is given by
$$
Y=
5\,+\,5\sum_{n\geq 1}\,b_n\,n^3\,\frac{q^n}{1-q^n}\,,
$$
$$
b_1 = 575\,,\;
b_2 = 121850\,,\;
b_3 = 63441275\,,\;
b_4 = 48493506000\,,\;
b_5 = 45861177777525\,.
$$
The prepotential $Z$ is defined (up to terms of order $\leq 2$ in $\log q$) by
$
Y=\left(q\frac{d}{dq}\right)^3 Z
$
and has an expansion
$$
Z(q)=\frac{5}{6}\log^3 q+
5\sum_{n\geq 1}b_n \mathrm{Li}_3(q^n)
$$
where $\mathrm{Li}_3(x)=\sum_{j\geq 1}\frac{x^j}{j^3}$ is the trilogarithm
function.
According to the Mirror Symmetry Conjecture \cite{C,M} the number $5b_n$ is the
number of rational curves of degree $n$ on a generic quintic hypersurface in
$\PP^4$. For the Integrality Conjecture one is interested in the less
astonishing, but still highly non-trivial fact that these numbers are integers.
Now one may note that the numbers $5b_n$ are integers if and only if for every
prime number $p$
\begin{equation}\label{eq:exa.quintic}
Z(q)-p^{-3} Z(q^p)\in \ZZ_p[[q]]
\end{equation}
where $\ZZ_p$ is the ring of $p$-adic integers.

In this paper we show that $p$-adic analogues of the integrality conjectures,
in particular an analogue of (\ref{eq:exa.quintic}), hold true
for ordinary CY3 crystals; see Formulas
(\ref{eq:qprops}) and (\ref{eq:yukinteger})
\footnote{
For some examples it has been proved in \cite{LY} that the expansion
coefficients of the mirror map (i.e. the coordinate transformation to canonical
coordinates)
are really integers, in $\ZZ$.}.
Throughout this paper the prime $p$ is fixed and the constructions
involve choices which may depend on $p$. For the integrality conjectures for
families of Calabi-Yau threefolds, like the above quintics, one must however
deal with all primes. Thus one is confronted with the challenge  to `glue' or
to `synchronize' those choices for the various
primes $p$. In \cite{S} we will describe the beginnings of such a synchronized
set-up for ``ordinariness''.
The present paper starts at the fairly abstract heights of $F$-crystals
and wants to find out from there whether the path of ordinariness leads to
interesting vistas. In \cite{S} we will start at a less abstract more geometric
level.

\section{Ordinary Hodge $F$-crystals}
\label{section 1}
We recall the theory of ordinary Hodge $F$-crystals; the general reference  is
\cite{D1}. We work over a perfect field $k$ of characteristic $p>0$ and
explicitly keep track of which results require passing to the algebraic closure
$\ok$ of $k$.
Let $W$ be the ring of Witt vectors of $k$ and $K$ the field of fractions of
$W$. So $W$ is a local ring with residue field $k$;
its field of fractions $K$ has characterisitic $0$;
its maximal ideal is $pW$ and $W$ is complete and separated for the $p$-adic
topology:  $\displaystyle{W=\lim_{\leftarrow m}W/p^mW}$. Moreover we need the
rings of formal power series
\begin{equation}\label{eq:A}
A_0=k[[t_1,\ldots,t_n]]\:,\quad
A=W[[t_1,\ldots,t_n]] .
\end{equation}

A \emph{crystal over $A_0$} is, by definition, a finitely generated free
$A$-module $H$ together with a connection
$$
\nabla:\,H\longrightarrow \OM_{A/W}^1\otimes_A H
$$
which is integrable and $p$-adically topologically nilpotent; this means
that if we set
$
D_i=(\frac{d}{dt_i}\otimes 1)\circ\nabla:H\rightarrow H,
$
with the derivation $\frac{d}{dt_i}$ viewed as a linear map
$\OM_{A/W}^1\rightarrow A$, then $D_iD_j=D_jD_i$
and $\displaystyle{\lim_{m\to\infty}D_i^m=0}$
in the $p$-adic topology on $\mathrm{End}_W(H)$.

On the characteristic $p$ rings $k$ and $A_0$ one has the \emph{Frobenius
endomorphism} $\sigma$ raising elements to their $p$-th power: $\sigma(x)=x^p$.
On the perfect field $k$ this is an automorphism and
it lifts canonically to an automorphism on $W$, also denoted by $\sigma$. There
are many different \emph{lifts of Frobenius on $A$}, i.e. ring endomorphisms
$\phi$ of $A$ which restrict to $\sigma$ on $W$
and reduce modulo $p$ to $\sigma$ on $A_0$.
If $H$ is an $A$-module and $\phi:A\rightarrow A$ is a lift of Frobenius we
write
$$
\phi^*H = A_\phi\otimes_A H
$$
where $A_\phi$ is $A$ viewed as an $A$-$A$-bimodule with left structure
via the identity map $id:A\rightarrow A$ and right structure via
$\phi:A\rightarrow A$. The natural map
$$
\phi^*: H\rightarrow \phi^*H \,,\qquad
\phi^*(h)=1\otimes h
$$
is $\phi$-linear, i.e.
$
\phi^*(a_1h_1+a_2h_2)=\phi(a_1)\phi^*(h_1)+\phi(a_2)\phi^*(h_2)$.

Now let $(H,\nabla)$ be a crystal as above and let
$\phi,\psi:A\rightarrow A$ be two lifts of Frobenius. Since $\phi$ and $\psi$
are equal modulo $p$ the connection $\nabla$ provides, via a kind of Taylor
expansion, a canonical isomorphism of $A$-modules
$$
\chi(\phi,\psi): \phi^*H\stackrel{\simeq}{\longrightarrow} \psi^*H;
$$
more precisely, the map $\chi(\phi,\psi)\phi^*:H\rightarrow\psi^*H$ is given by
the formula
$$
\chi(\phi,\psi)\phi^*(x)=\sum_{m_1,\ldots,m_n\geq 0}
\prod_{j=1}^n\frac{p^{m_j}}{m_j!}
\left(\frac{\phi (t_j)-\psi(t_j)}{p}\right)^{m_j}\psi^*(D_1^{m_1}
\circ\ldots\circ D_n^{m_n}(x)).
$$
Note that for $m\geq 1$ the rational number $\frac{p^m}{m!}$ has $p$-adic
valuation $> m-\frac{m}{p-1}\:\geq 0$.

\begin{definition}\label{Fcrystal}
One says that a crystal $(H,\nabla)$ is an \emph{$F$-crystal over $A_0$} if for
every
lift of Frobenius $\phi:A\rightarrow A$ there is given a homomorphism of
$A$-modules
\begin{equation}\label{eq:frobcrys}
F(\phi):\,\phi^* H\rightarrow H
\end{equation}
which is horizontal for the connection $\nabla$, i.e. the square
\begin{equation}\label{eq:horizontal}
\begin{array}{rcccl}
&H&\stackrel{\nabla}{\longrightarrow}&
\OM_{A/W}^1\otimes H&\\[.6ex]
\mbox{\scriptsize{$F(\phi)\phi^*$}}\hspace{-1em}&\downarrow&&\downarrow&
\hspace{-2.5em} \mbox{\scriptsize{$\phi\otimes F(\phi)\phi^*$}} \\[.6ex]
&H&\stackrel{\nabla}{\longrightarrow}&\OM_{A/W}^1\otimes H&
\end{array}
\end{equation}
is commutative, and such that for every pair of lifts of Frobenius
$\phi,\psi:A\rightarrow A$
\begin{equation}\label{eq:Fvar}
F(\psi)\circ\chi(\phi,\psi)=F(\phi).
\end{equation}
Moreover $F(\phi)\otimes \QQ_p:\,\phi^* H\otimes \QQ_p
\rightarrow H\otimes \QQ_p$ should be an isomorphism.
\\
If for one, and hence every, lift of Frobenius
$\phi:A\rightarrow A$ the homomorphism $F(\phi)$ in (\ref{eq:frobcrys}) is an
isomorphism one says that $H$ is a \emph{unit $F$-crystal}.
\end{definition}
Combining Formula (\ref{eq:Fvar}) with the Taylor expansion formula for
$\chi(\phi,\psi)$ one gets
\begin{equation}\label{eq:FTaylor}
F(\phi)\phi^*(x)=\sum_{m_1,\ldots,m_n\geq 0}
\prod_{j=1}^n\frac{p^{m_j}}{m_j!}
\left(\frac{\phi (t_j)-\psi(t_j)}{p}\right)^{m_j}
F(\psi)\psi^*(D_1^{m_1}
\circ\ldots\circ D_n^{m_n}(x)).
\end{equation}

The horizontality relation for $F(\phi)$ and $\nabla$ (\ref{eq:horizontal})
makes it possible to solve the problem of finding $\nabla$-flat sections
(i.e. essentially solving differential equations) by
finding $F(\phi)\phi^*$-fixed sections (i.e. essentially solving polynomial
equations). Indeed, one has the following result.

\begin{lemma}\label{eq:solve}
Suppose $e\in H$ is such that $F(\phi)\phi^*e=e$ for some lift of Frobenius
$\phi:A\rightarrow A$. Then $\nabla e =0$ and $F(\psi)\psi^*e=e$ holds for
every lift of Frobenius $\psi$.
\end{lemma}
\begin{proof}
By definition there is for every $a\in A$ an element $a_1\in A$ such that
$\phi(a)=a^p+pa_1$. From this one gets by induction for every $a\in A$ and
$m\geq 1$ elements $a_1,\ldots,a_m\in A$ such that
$$
\phi^m(a)=a^{p^m}+p a_1^{p^{m-1}}+\ldots+p^ma_m
$$
and hence
$
\phi^m(bda):=\phi^m(b)d(\phi^m(a))\in p^m\OM_{A/W}^1
$
for every $1$-form $bda$.\\
If $F(\phi)\phi^*e=e$, then (\ref{eq:horizontal}) yields for every $m\geq 1$
$$
\nabla e=(\phi^m\otimes (F(\phi)\phi^*)^m)\circ\nabla e\in
p^m\OM_{A/W}^1\otimes H.
$$
Taking $m\rightarrow\infty$ we see $\nabla e=0$.
The equality $F(\psi)\psi^*e=e$ for any lift of Frobenius $\psi$ now follows
from (\ref{eq:FTaylor}).
\end{proof}

Let $H$ be an $F$-crystal\footnote{For simplicity, we do not
explicitly mention the connection $\nabla$ and the maps $F(\phi)$.}
and $H_0=H\otimes_A A_0$. For a lift of Frobenius $\phi:A\rightarrow A$ and an
integer $i$ we set
\begin{eqnarray*}
\Fil_iH_0&=&\{x\in H_0\,|\,\exists y\in H\textrm{ lifting } x
\textrm{ s.t. } p^iy\in\mathrm{Im} F(\phi)\},
\\
\Fil^iH_0&=&\{x\in H_0\,|\,\exists y\in H\textrm{ lifting } x
\textrm{ s.t. } F(\phi)\phi^*y\in p^iH \}.
\end{eqnarray*}
These are $A_0$-submodules of $H_0$. One can show that they are independent of
the chosen lift of Frobenius $\phi$; see \cite{D1} \S 1.3.
Trivially, $\Fil_iH_0=0$ for $i\leq -1$ and
$\Fil^iH_0=H_0$ for $i\leq 0$. Moreover, it can be shown that there is an
integer $N$ such that $\Fil_NH_0=H_0$ and  $\Fil^{N+1}H_0=0$.
One then says that \emph{$H$ is of level $\leq N$}.
Thus, one finds the filtrations
\begin{eqnarray}
\label{eq:filcon2}
0=\Fil_{-1}H_0\subset\Fil_0H_0\subset\Fil_1H_0\subset&\ldots&\subset
\Fil_{N-1}H_0\subset\Fil_NH_0=H_0
\\
\label{eq:filhodge2}
H_0=\Fil^0H_0\supset\Fil^1H_0\supset\Fil^2H_0\supset&\ldots&
\supset\Fil^NH_0\supset\Fil^{N+1}H_0=0
\end{eqnarray}
(\ref{eq:filcon2}) is called the \emph{conjugate filtration} and
(\ref{eq:filhodge2}) is called the \emph{Hodge filtration}.
The connection $\nabla$ on $H$ induces a connection $\nabla_0$ on $H_0$. The
conjugate filtration is horizontal for $\nabla_0$ and the Hodge filtration
satisfies the \emph{Griffiths transversality condition}, i.e. for all $i$:
$$
\nabla_0\Fil_i H_0\subset\OM_{A_0/k}^1\otimes\Fil_i H_0
\qquad\textrm{and}\qquad
\nabla_0\Fil^i H_0\subset \OM_{A_0/k}^1\otimes\Fil^{i-1} H_0
{}.
$$
Let $\gr_\bullet H_0$ and $\gr^\bullet H_0$ denote the graded modules
associated with the conjugate and Hodge filtrations respectively. If
$\gr_\bullet H_0$ is a free $A_0$-module, then so is $\gr^\bullet H_0$ and for
every $i$ the modules
$\gr_iH_0=\Fil_iH_0/\Fil_{i-1}H_0$ and $\gr^iH_0=\Fil^iH_0/\Fil^{i+1}H_0$ are
canonically isomorphic. The rank of the free $A_0$-module
$\gr_iH_0\simeq\gr^iH_0$ is called \emph{the $i$-th Hodge number of $H$}.

\begin{definition}\label{defHodgecrystal}
A \emph{Hodge $F$-crystal} over $A_0$ is an $F$-crystal $H$ over $A_0$
equipped with a filtration by free $A$-submodules
\begin{equation}\label{eq:Hodgelift}
H=\Fil^0H\supset\Fil^1H\supset\ldots
\supset\Fil^iH\supset\Fil^{i+1}H\supset\ldots
\end{equation}
(called the Hodge filtration on $H$) which lifts the Hodge filtration from
$H_0$ and satisfies Griffiths transversality, i.e. for every $i$
\begin{equation}\label{eq:transversality}
\Fil^iH\otimes_A A_0=\Fil^iH_0
\qquad\textrm{and}\qquad
\nabla\Fil^i H\subset \OM_{A/W}^1\otimes_A \Fil^{i-1} H.
\end{equation}
\end{definition}

\begin{definition}\label{defordinary}
An $F$-crystal $H$ over $A_0$ is said to be \emph{ordinary} if
$\gr_\bullet H_0$ is a free $A_0$-module and the conjugate filtration and Hodge
filtration on $H_0$ are opposite, i.e.
$$
H_0=\Fil_iH_0\oplus\Fil^{i+1}H_0 \qquad\textrm{ for every } i.
$$
\end{definition}

\begin{proposition}\label{propordinary}
\textup{(\cite{D1} prop.1.3.2)}
Let $H$ be an $F$-crystal over $A_0$ such that $\gr_\bullet H_0$ is a free
$A_0$-module. Then $H$ is ordinary if and only if there is a filtration by
sub-$F$-crystals
\begin{equation}\label{eq:Ufilt}
0=U_{-1}\subset U_0\subset U_1\subset\ldots\subset U_i\subset
U_{i+1}\subset\ldots
\end{equation}
such that for every $i$
$$
U_i\otimes_A A_0=\Fil_iH_0\quad\textrm{ and }\quad
U_i/U_{i-1}\simeq V_i(-i)
$$
where $V_i$ is a unit $F$-crystal and $(-i)$ is Tate twist, i.e.
$V_i(-i)$ is the same $A$-module with connection as $V_i$, but for every lift
of Frobenius $\phi$ the map $F(\phi)\phi^*$ on $V_i(-i)$ is $p^iF(\phi)\phi^*$
on $V_i$.
The filtration (\ref{eq:Ufilt}) is unique.
\qed
\end{proposition}
\begin{proposition}\label{propHodgedeco}
\textup{(\cite{D1} prop.1.3.6)}
For an ordinary Hodge $F$-crystal $H$ the filtrations $U_\bullet$ and
$\Fil^\bullet H$ are opposite, i.e. for every $i$
$$
H=U_i\oplus \Fil^{i+1}H.
$$
As a consequence one has a decomposition
$$
H=\oplus_i H^i\;,\qquad H^i=U_i\cap\Fil^iH.
$$
\qed
\end{proposition}

The above proposition gives a first result on the structure of
ordinary Hodge $F$-crystals.
Stronger results (Theorem \ref{thmstructure} below) can be obtained over the
algebraic closure $\ok$ of the base field $k$.
Let $\oW$ the ring of Witt vectors of $\ok$ and
$$
\oA_0=\ok[[t_1,\ldots,t_n]]\,,\qquad\oA=\oW[[t_1,\ldots,t_n]]
\,,\qquad \Ku=\oW[\textstyle{\frac{1}{p}}].
$$
One can base change an $F$-crystal $H$ over $A_0$ to an $F$-crystal $\oH$ over
$\oA_0$ as follows. By tensoring with $\oA$ one gets
the free $\oA$-module with connection
$$
\oH=\oA\otimes_A H\,,\qquad\nabla: \oH\rightarrow \OM_{\oA/\oW}^1\otimes_\oA
\oH.
$$
If $\phi:\oA\rightarrow \oA$ is a lift of Frobenius such that
$\phi(A)\subset A$, then the $A$-linear map $F(\phi):\phi^*H\rightarrow H$
gives the $\oA$-linear map
$F(\phi):\phi^*\oH\rightarrow \oH$.
If $\psi:\oA\rightarrow \oA$ is an arbitrary lift of Frobenius,
the connection $\nabla$ provides an isomorphism
$\chi(\psi,\phi):\psi^*\oH\stackrel{\simeq}{\longrightarrow}\phi^*\oH$
and we can define
$F(\psi)=F(\phi)\circ\chi(\psi,\phi):\psi^*\oH\rightarrow\oH$.

\begin{definition}\label{defdivisible}
Let $\oH$ be a Hodge $F$-crystal over $\oA_0$ with Hodge filtration
$\{\Fil^i\oH\}_{0\leq i\leq N}$.
One says that $\oH$ is \emph{divisible} if for some lift of Frobenius $\phi$
\begin{equation}\label{eq:divisible}
F(\phi)\phi^* \Fil^i\oH\subset  p^i\oH\qquad\textrm{for}\quad i=0,1,\ldots,N.
\end{equation}
If $p>N$ and (\ref{eq:divisible}) holds for one lift of Frobenius,
then it holds for every lift of Frobenius (see \cite{K1} \S 5.0).
\end{definition}

Th\'eor\`eme 1.4.2 of \cite{D1} describes the structure of ordinary Hodge
$F$-crystals of level $\leq 1$.
The following theorem generalizes that result to higher levels.

\begin{theorem}\label{thmstructure}
Let $\oH$ be a divisible ordinary Hodge $F$-crystal over $\oA_0$ of level $\leq
N<p$ with filtrations $U_\bullet$ and $\Fil^\bullet H$
as in (\ref{eq:Ufilt}) and (\ref{eq:Hodgelift}) respectively.
Then
$$
\oH=\bigoplus_{i=0}^N \oH^i\:,\qquad \oH^i=\oU_i\cap\Fil^i\oH
$$
and
there is a basis $\{e_{im}\}_{0\leq i\leq N, 1\leq m\leq h^i}$
of $\oH$ with $e_{im}\in\oH^i$ for all $i,m$ and there is a matrix $\cT$ with
entries in $\Ku[[t_1,\ldots,t_n]]$, such that
\begin{itemize}
\item
$\cT=(\tau_{ij})_{0\leq i,j\leq N}$ with
$\tau_{ij}$ an $h^i\times h^j$-matrix,
$\tau_{ij}=0$ for $j<i$,\\
$\tau_{ii}=$identity-matrix,
the matrix $\tau_{ij}(0)$ has entries in
$p^{j-i}\oW$.
\item
the matrix of the connection $\nabla$ with respect to the basis
$\{e_{im}\}$ is
\begin{equation}\label{eq:conmat}
\cT^{-1}\cdot d\cT
\end{equation}
\item
for every lift of Frobenius $\psi$ the matrix of the map $F(\psi)\psi^*$ with
respect to the basis
$\{e_{im}\}$ is
\begin{equation}\label{eq:frobmat}
\cT^{-1}\cdot P\cdot\psi(\cT)
\end{equation}
where $P$ is the diagonal matrix
\begin{equation}\label{eq:P}
P=\mathrm{diag}(1,\ldots,1,p,\ldots,p,\ldots,p^j,\ldots,
p^j,\ldots,p^N,\ldots,p^N)
\end{equation}
with the entry $p^j$ repeated $h^j$-times.
\end{itemize}
\end{theorem}
\begin{proof}
For $i=0,\ldots,N$ there is an isomorphism of $\oA$-modules
$\oH^i\simeq \oU_i/\oU_{i-1}=\oV_i(-i)$ with $\oV_i$ a unit $F$-crystal.
According to \cite{D1} prop.1.2.2 there is an $\oA$-basis
$f_{i1},\ldots,f_{ih^i}$ of $\oV_i$ such that
\begin{equation}\label{eq:unitbase}
\nabla f_m=0\;,\qquad F(\psi)\psi^*f_m=f_m\qquad\textrm{for}\quad
m=1,\ldots,h^i
\end{equation}
for every lift of Frobenius $\psi:\oA\rightarrow\oA$.
Lift the basis $f_{i1},\ldots,f_{ih^i}$ of $\oV_i$ to a basis
$e_{i1},\ldots,e_{ih^i}$ of $\oH^i$.
This gives the basis $\{e_{im}\}_{0\leq i\leq N, 1\leq m\leq h^i}$ for $\oH$.

Equations
(\ref{eq:unitbase}) and Griffiths transversality show that the connection
matrix with respect to this basis has the following shape
\begin{equation}\label{eq:matrix1}
M_\nabla=\left(
\begin{array}{cccccc}
0&\eta_{01}&0&0&\ldots&0\\0&0&\eta_{12}&0&\ldots&0\\
&&\ddots&\ddots&&
\\ 0&\ldots&\ldots&0&\eta_{N-2,N-1}&0\\
0&\ldots&\ldots&\ldots&0&\eta_{N-1,N}\\
0&\ldots&\ldots&\ldots&0&0
\end{array}
\right)
\end{equation}
with $\eta_{i,i+1}$ a matrix of size $h^i\times h^{i+1}$ with entries in
$\OM_{\oA/\oW}^1$. Integrability of the connection implies for all $i$
\begin{equation}\label{eq:integrability}
d\eta_{i,i+1}=0\qquad\textrm{and}\qquad\eta_{i,i+1}\eta_{i+1,i+2}=0.
\end{equation}
 The entries of
$\eta_{i,i+1}$ are closed $1$-forms and hence exact, by the Poincar\'e lemma.
So
there is a matrix $\tau_{i,i+1}$ of size $h^i\times h^{i+1}$ with entries in
$\Ku[[t_1,\ldots,t_n]]$ such that
\begin{equation}\label{eq:etatau1}
\eta_{i,i+1}=d\tau_{i,i+1}.
\end{equation}
Using the equations in (\ref{eq:integrability}) and induction on $j-i$ one sees
that for $i+2\leq j\leq N$ there exist matrices $\tau_{ij}$ of size $h^i\times
h^j$ with entries in $\Ku[[t_1,\ldots,t_n]]$ such that
\begin{equation}\label{eq:etatau2}
\tau_{i,j-1}\eta_{j-1,j}=d\tau_{ij}.
\end{equation}
For $j<i$ we set $\tau_{ij}=h^i\times h^j$-zero-matrix
and for $j=i$ we set $\tau_{ii}=h^i\times h^i$-identity-matrix.
We collect the matrices $\tau_{ij}$ into one (block structured) matrix
$$
\cT=(\tau_{ij})_{0\leq i,j\leq N}.
$$
Then (\ref{eq:etatau1}) and (\ref{eq:etatau2}) are equivalent with
\begin{equation}\label{eq:dT}
\cT M_\nabla =d\cT.
\end{equation}
We fix the constants of integration by taking
the constant term $\cT (0)$ of $\cT$ as in (\ref{eq:T0}) below.
Take the lift of Frobenius $\phi:\oA\rightarrow\oA$ given by
\begin{equation}\label{eq:lift}
\phi(\sum a_{m_1,\ldots,m_n}t_1^{m_1}\cdots t_n^{m_n})=
\sum \sigma(a_{m_1,\ldots,m_n})t_1^{pm_1}\cdots t_n^{pm_n}
\end{equation}
with $\sigma$ the standard Frobenius map on $\oW$.
Equations (\ref{eq:unitbase}) and the assumption on divisibility show
$$
F(\phi)\phi^*e_{im}-p^ie_{im}\in p^i\oU_{i-1}.
$$
Therefore the matrix $M_\phi$ of the $\phi$-linear map
$F(\phi)\phi^*:\oH\rightarrow\oH$ has the shape
$$
M_\phi=LP
$$
with $P$ as in (\ref{eq:P})
and $L$ an uppertriangular matrix with entries in $\oA$ and $1$'s along the
diagonal. Let $L(0)$ denote the constant term of $L$.
After these preparations we define $\cT$ to be the unique matrix which
satisfies the differential equation (\ref{eq:dT}) and the initial condition
\begin{equation}\label{eq:T0}
\cT (0)=\prod_{m=1}^\infty \sigma^{-m}(P^{-m}L(0)P^{m});
\end{equation}
this product (with the factors ordered from right to left for increasing $m$)
converges since the matrix $P^{-m}L(0)P^{m}-I$ has entries in $p^{m}\oW$.

The fact that $F(\phi)$ is horizontal for the connection $\nabla$
(see diagram (\ref{eq:horizontal})) is expressed by the matrix equation
$
M_\phi\phi(M_\nabla) = M_\nabla M_\phi+d M_\phi.
$
Multiplying this matrix equation by $\cT$ and using Equation (\ref{eq:dT})
gives
$$
\cT M_\phi \phi(\cT)^{-1}d\phi(\cT)=d( P^{-1} \cT M_\phi)
$$
and hence $\cT M_\phi \phi(\cT)^{-1}$ is constant.
Since $\phi(\cT)(0)=P^{-1}\cT(0) M_\phi(0)$
by (\ref{eq:T0}) we conclude:
$$
M_\phi=\cT ^{-1}\,P\,\phi(\cT).
$$
This shows that (\ref{eq:frobmat}) holds for the particular lift of Frobenius
$\phi$. In order to check it for an arbitrary lift of Frobenius we first define
the new basis
$\{\ep_{im}\}_{0\leq i\leq N, 1\leq m\leq h^i}$
for $\oH\otimes\Ku[[t_1,\ldots,t_n]]$ so that $\cT^{-1}$ is the matrix whose
columns give the coordinates of the new basis vectors with respect to the old
basis
$\{e_{im}\}_{0\leq i\leq N, 1\leq m\leq h^i}$:
\begin{equation}\label{eq:newbasetrafo}
\{\ep_{im}\}=\{e_{im}\}\cT^{-1}.
\end{equation}
The triangular shape of the matrix $\cT$ implies
\begin{equation}\label{eq:filnewbase}
\ep_{im}\in \oU^i\otimes\Ku[[t_1,\ldots,t_n]]
\qquad\textrm{for all}\quad i,m.
\end{equation}
The connection matrix with respect to the new basis is
$d(\cT^{-1})+M_\nabla\cT^{-1}=0$. So
\begin{equation}\label{eq:dnewbase}
\nabla\ep_{im}=0\qquad\textrm{for all}\quad i,m.
\end{equation}
The matrix of $F(\phi)\phi^*$ with respect to the new basis is
$\cT M_\phi\phi(\cT)^{-1}=P$. So
$
F(\phi)\phi^*\ep_{im}=p^i\ep_{im}$ for all $i,m$.
This holds for the special lift of Frobenius $\phi$ (\ref{eq:lift}).
If $\psi$ is any lift of Frobenius (\ref{eq:FTaylor})
and (\ref{eq:dnewbase})
show $F(\psi)\psi^*\ep_{im}=F(\phi)\phi^*\ep_{im}$. Thus
for any lift of Frobenius $\psi$ we have
\begin{equation}\label{eq:Fnewbase}
F(\psi)\psi^*\ep_{im}=p^i\ep_{im}\qquad\textrm{for all}\quad i,m.
\end{equation}
The matrix $M_\psi$ of $F(\psi)\psi^*$ with respect to the basis $\{e_{im}\}$
is therefore
$$
M_\psi=\cT^{-1}\cdot P\cdot\psi(\cT),
$$
as claimed in (\ref{eq:frobmat}).
To finish the proof of Theorem \ref{thmstructure} we note that (\ref{eq:T0})
implies that the matrix $\tau_{ij}(0)$ has entries in
$p^{j-i}\oW$ for $j>i$.
\end{proof}
\section{Ordinary CY3-crystals over an algebraically closed
field}\label{sectionCY3}
We now specialize the general result of Theorem \ref{cy3structure} to crystals
of the type that appears as the third crystalline cohomology of an ordinary
Calabi-Yau threefold over the algebraically closed field $\ok$.
More precisely, we want to see what happens, if we combine the hypotheses in
Theorem \ref{cy3structure} with the hypotheses for the variation of Hodge
structure of Calabi-Yau threefolds, which are explicitly stated in \cite{BG} \S
2. So \emph{from now on
$p>3$ and $\oH$ is a divisible ordinary Hodge $F$-crystal over $\oA_0$ of level
$3$} with filtrations
\begin{eqnarray*}
&&0\subset \oU_0\subset \oU_1\subset \oU_2\subset \oU_3=\oH\\
&&\oH=\Fil^0\oH\supset\Fil^1\oH\supset\Fil^2\oH\supset
\Fil^3\oH\supset 0.
\end{eqnarray*}
such that in the decomposition
$\oH=\oH^0\oplus\oH^1\oplus\oH^2\oplus\oH^3\,,\quad \oH^i=\oU_i\cap\Fil^i\oH$,
\begin{equation}\label{eq:Hodgenumbers}
\rank\,\oH^0=\rank\,\oH^3=1\qquad\textrm{and}\qquad\rank\,\oH^1=\rank\,\oH^2=h;
\end{equation}
$\oH^i$ would be denoted as $H^{i,3-i}$ in the Hodge theoretic setting of
\cite{BG}.
\\
Moreover we assume that there is given a non-degenerate alternating bilinear
form
$$
\langle\:,\:\rangle\,:\:\oH\times\oH\rightarrow\oA
$$
such that for all $x,y\in\oH$ and for every lift of Frobenius
$\phi:\oA\rightarrow\oA$
\begin{equation}\label{eq:formconditions}
\begin{array}{rcl}
\langle\nabla x,y\rangle + \langle x,\nabla y\rangle &=& d\langle x,y\rangle
\\[.5ex]
\langle F(\phi)\phi^* (x),F(\phi)\phi^*(y)\rangle &=& p^3\phi(\langle
x,y\rangle).
\end{array}
\end{equation}
We also require the Riemann bilinear relations:
\begin{equation}\label{eq:Riemann}
(\Fil^3\oH)^\perp=\Fil^1\oH\:,\qquad (\Fil^2\oH)^\perp=\Fil^2\oH\:.
\end{equation}
\begin{definition}\label{def CY3crys}
An \emph{ordinary CY3 crystal} over $\oA_0$ is a divisible ordinary Hodge
$F$-crystal of level $3$
over $\oA_0$ equipped with a non-degenerate alternating bilinear form
$\langle\:,\:\rangle$
such that (\ref{eq:Hodgenumbers}), (\ref{eq:formconditions}),
(\ref{eq:Riemann}) are satisfied.
\end{definition}

In view of the Griffiths tranversality condition (\ref{eq:transversality}) the
connection $\nabla$ induces an $\oA$-linear map
\begin{equation}\label{eq:conngrhodge}
T_{\oA/\oW}\rightarrow \mathrm{Hom}(\oH^3,\oH^2)
\end{equation}
where $T_{\oA/\oW}=\mathrm{Hom}(\OM_{\oA/\oW}^1,\oA)$.

\begin{theorem}\label{cy3structure}
Let $p>3$ and let $\oH$ be an ordinary CY3 crystal over $\oA_0$.
Then there is a basis $\{e_i\}_{i=0,\ldots,2h+1}$, such that
$$
e_0\in\oH^0,\qquad
e_1,\ldots,e_h\in\oH^1,\qquad
e_{h+1},\ldots,e_{2h}\in\oH^2,\qquad e_{2h+1}\in\oH^3.
$$
and such that the Gramm matrix of the form $\langle\:,\:\rangle$ with respect
to this basis is
\begin{equation}\label{eq:eGramm}
(\langle e_i,e_j\rangle)_{0\leq i,j\leq 2h+1}
=
\left(\begin{array}{cccc}
0&0&0&-1\\0&0&I_h&0\\0&-I_h&0&0\\ 1&0&0&0
\end{array}\right)
\end{equation}
with $I_h= h\times h$-identity-matrix.
Moreover there is a matrix $\cT$ such that
the matrix of the connection $\nabla$ with respect to the basis
$e_0,\ldots,e_{2h+1}$ is
\begin{equation}\label{eq:conmat2}
\cT^{-1}\cdot d\cT
\end{equation}
and such that for every lift of Frobenius $\psi$ the matrix of the map
$F(\psi)\psi^*$ with respect to the basis
$e_0,\ldots,e_{2h+1}$ is
\begin{equation}\label{eq:frobmat2}
\cT^{-1}\cdot P\cdot\psi(\cT)
\end{equation}
with $P$ the (block) diagonal matrix $\textrm{diag}(1,pI_h,p^2I_h,p^3)$.

This matrix $\cT$
factors as \footnote{$\qquad.^\star$ means matrix transpose}
\begin{equation}\label{eq:Tfac}
\cT
=
\left(\begin{array}{cccc}
1&\tau_{23}^\star&0&0\\0&I_h&0&0
\\ 0&0&I_h&\tau_{23}\\0&0&0&1
\end{array}\right)
\left(\begin{array}{cccc}
1&0&-\tau_{13}^\star&Z\\0&I_h&\tau_{12}&\tau_{13}
\\ 0&0&I_h&0\\0&0&0&1
\end{array}\right)
\end{equation}
where $Z,\,\tau_{23},\,\tau_{13},\,\tau_{12}$ are matrices with entries
in $\Ku[[t_1,\ldots,t_n]]$ of respective sizes $1\times 1$, $h\times 1$,
$h\times 1$, $h\times h$.

Assume in  addition to the preceding assumptions that the map in
(\ref{eq:conngrhodge}) is an isomorphism, and hence  $n=h$.
Let
$\tau_1,\ldots,\tau_h\in\Ku[[t_1,\ldots,t_h]] $ denote the components of the
$h\times 1$-matrix $\tau_{23}$:
\begin{equation}\label{eq:taudef}
\tau_{23}=(\tau_1,\ldots,\tau_h)^\star.
\end{equation}
Then the matrix
\begin{equation}\label{eq:dtaudt}
\left(\frac{\partial\tau_i}{\partial t_j}\right)_{1\leq i,j\leq h}
\end{equation}
is invertible over $\Ku[[t_1,\ldots,t_h]] $ and
\begin{eqnarray}
\label{eq:tau13b}
\tau_{13}&=&-\frac{1}{2}\left(\frac{\partial Z}{\partial \tau_1},
\ldots,\frac{\partial Z}{\partial \tau_h}\right)^\star
\\
\label{eq:tau12b}
\tau_{12}&=&-\frac{1}{2}\left(
\frac{\partial^2 Z}{\partial \tau_i\partial \tau_j}
\right)_{1\leq i,j\leq h}.
\end{eqnarray}
\end{theorem}
\begin{proof}
Take any $\oA$-basis of $\oH$ and matrix for which the statements in Theorem
\ref{thmstructure} hold and denote these as
$e'_0,e'_1,\ldots,e'_h,e'_{h+1},\ldots,e'_{2h},e'_{2h+1}$ and $\cT'$
respectively. So
$$
e'_0\in\oH^0,\qquad
e'_1,\ldots,e'_h\in\oH^1,\qquad
e'_{h+1},\ldots,e'_{2h}\in\oH^2,\qquad e'_{2h+1}\in\oH^3.
$$
Also take the basis
$\ep'_0,\ep'_1,\ldots,\ep'_h,\ep'_{h+1},\ldots,\ep'_{2h},\ep'_{2h+1}$ as in
(\ref{eq:filnewbase}), related to the previous basis by the transformation with
the matrix $\cT'$; see (\ref{eq:newbasetrafo}).

Then $\nabla\ep'_i=\nabla\ep'_j=0$ and hence
$d\langle\ep'_i,\ep'_j\rangle=\langle\nabla\ep'_i,\ep'_j\rangle+
\langle\ep'_i,\nabla\ep'_j\rangle=0$ for all $i,j$.
This implies
$\langle\ep'_i,\ep'_j\rangle\in\oW$ for all $i,j$.
Moreover we have $F(\phi)\phi^* (\ep'_i)=p^r \ep'_i$ with
$r=\lceil\!\frac{i}{h}
\!\rceil:=$ smallest integer $\geq\frac{i}{h}$,
and
$\langle F(\phi)\phi^* (\ep'_i),F(\phi)\phi^* (\ep'_j)\rangle
=p^3\sigma(\langle\ep'_i,\ep'_j\rangle)$.
This implies
$\sigma(\langle\ep'_i,\ep'_j\rangle)=\langle\ep'_i,\ep'_j\rangle$
and hence $\langle\ep'_i,\ep'_j\rangle\in\ZZ_p$
if $\lceil\!\frac{i}{h}\!\rceil+\lceil\!\frac{j}{h}\!\rceil= 3$,
whereas
$\langle\ep'_i,\ep'_j\rangle=0$ if $\lceil\!\frac{i}{h}\!\rceil+
\lceil\!\frac{j}{h}\!\rceil\neq 3$.
Now take the (block structured) matrix with entries in $\ZZ_p$:
$$
M=\left(\begin{array}{cccc}
1&0&0&0\\0&
I_h&0&0\\0&0&
(\langle\ep'_i,\ep'_{j+h}\rangle)_{1\leq i,j\leq h}
&0\\ 0&0&0&-\langle\ep'_0,\ep'_{2h+1}\rangle
\end{array}\right)
$$
This matrix is invertible because the form $\langle\:,\:\rangle$ is
non-degenerate. We use it for the basis transformations
\begin{eqnarray*}
(e_0,e_1,\ldots,e_h,e_{h+1},\ldots,e_{2h},e_{2h+1})
&:=&(e'_0,e'_1,\ldots,e'_h,e'_{h+1},\ldots,e'_{2h},e'_{2h+1})M^{-1}\\
(\ep_0,\ep_1,\ldots,\ep_h,\ep_{h+1},\ldots,\ep_{2h},\ep_{2h+1})
&:=&(\ep'_0,\ep'_1,\ldots,\ep'_h,\ep'_{h+1},\ldots,\ep'_{2h},\ep'_{2h+1})
M^{-1}
\end{eqnarray*}
For the new basis
$\ep_0,\ldots,\ep_{2h+1}$, the (block structured) Gramm matrix is
\begin{equation}\label{eq:goodgramm}
(\langle\ep_i,\ep_j\rangle)_{0\leq i,j\leq 2h+1}
=
\left(\begin{array}{cccc}
0&0&0&-1\\0&0&I_h&0\\0&-I_h&0&0\\ 1&0&0&0
\end{array}\right).
\end{equation}
We set
$$
\cT=M\cT'M^{-1},
$$
so that we have the basis relation
$$
(e_0,e_1,\ldots,e_h,e_{h+1},\ldots,e_{2h},e_{2h+1})=
(\ep_0,\ep_1,\ldots,\ep_h,\ep_{h+1},\ldots,\ep_{2h},\ep_{2h+1})\cT.
$$
The block structure of the matrix $\cT$ is
\begin{equation}\label{eq:goodT}
\cT
=
\left(\begin{array}{cccc}
1&\tau_{01}&\tau_{02}&\tau_{03}\\0&I_h&\tau_{12}&\tau_{13}
\\ 0&0&I_h&\tau_{23}\\0&0&0&1
\end{array}\right).
\end{equation}
The statements contained in formulas
(\ref{eq:dnewbase}), (\ref{eq:Fnewbase}), (\ref{eq:conmat}), (\ref{eq:frobmat})
are still valid with respect to the new bases.

The Gramm matrix of the form $\langle\,,\,\rangle$ with respect to the basis
$e_0,\ldots,e_{2h+1}$ is
\begin{eqnarray}\nonumber
&&\hspace{-6em}
(\langle e_i,e_j\rangle)_{0\leq i,j\leq 2h+1}
=
\\
\label{eq:gramme}
&=&
\left(\begin{array}{rrrr}
0&0&0&-1\\ 0&0&I_h&\tau_{23}-\tau_{01}^\star\\
0&-I_h&-\tau_{12}+\tau_{12}^\star&
-\tau_{13}+\tau_{12}^\star\tau_{23}-\tau_{02}^\star\\
1&
\tau_{01}-\tau_{23}^\star&\tau_{02}-\tau_{23}^\star\tau_{12}+
\tau_{13}^\star&
0
\end{array}\right)
\end{eqnarray}
The Riemann bilinear relations (\ref{eq:Riemann}) boil down to
$$
\langle e_i,e_{2h+1}\rangle=0\quad\textrm{for }\;i\geq 1\,,\qquad
\langle e_i,e_j\rangle=0\quad\textrm{for }\;i,j\geq h+1.
$$
In view of (\ref{eq:gramme}) this means firstly that the Gramm matrix is indeed
as in (\ref{eq:eGramm}) and secondly that
\begin{equation}\label{eq:taurel}
\tau_{23}=\tau_{01}^\star
\,,\qquad
\tau_{12}=\tau_{12}^\star
\,,\qquad
\tau_{02}=\tau_{23}^\star\tau_{12}-\tau_{13}^\star.
\end{equation}
The latter equalities are equivalent with the factorization of $\cT$ as in
(\ref{eq:Tfac}), with
\begin{equation}\label{eq:yuk}
Z=\tau_{03}-\tau_{23}^\star\tau_{13}.
\end{equation}
{}From Equations (\ref{eq:etatau1}), (\ref{eq:etatau2}),
(\ref{eq:taurel}) one gets immediately
\begin{equation}\label{eq:tau12}
d\tau_{13}= \tau_{12}d\tau_{23}
\end{equation}
and with some straightforward extra calculation
\begin{equation}\label{eq:tau13}
dZ=d(\tau_{03}-\tau_{13}^\star\tau_{23})=-2\tau_{13}^\star d\tau_{23}.
\end{equation}

Now assume that the map in (\ref{eq:conngrhodge}) is an isomorphism.
The partial derivations $\frac{\partial}{\partial t_j}$  ($j=1,\ldots,h$)
constitute a natural basis for $T_{\oA/\oW}$.
A natural basis for $\mathrm{Hom}(\oH^3,\oH^2)$ consists of the maps
$e_{2h+1}\mapsto e_{h+i}$ ($i=1,\ldots,h$).
With respect to these bases the map (\ref{eq:conngrhodge}) is given by
the matrix
\begin{equation}\label{eq:kodaira-spencer}
\left(\frac{\partial\tau_i}{\partial t_j}\right)_{1\leq i,j\leq h}.
\end{equation}
This matrix is therefore invertible.
Equation (\ref{eq:tau13b}) is an obvious restatement of (\ref{eq:tau13}), while
Equation (\ref{eq:tau12b}) follows from
(\ref{eq:tau12}) and (\ref{eq:tau13b}).
\end{proof}

In the preceding proof we have essentially followed the same algebraic path as
Bryant and Griffiths in their analysis of the variation of the Hodge structure
of (complex) Calabi-Yau threefolds \cite{BG,F}.
Formulas (\ref{eq:yuk}), (\ref{eq:tau13b}), (\ref{eq:tau12b}) are familiar from
the theory of Bryant and Griffiths. The following computation shows that the
equally familiar formula for the \emph{cubic form of Bryant and Griffiths}
(also known as \emph{Yukawa coupling}) holds in the present situation as well:

\begin{corollary}\label{coryukawa}
In the situation of Theorem \ref{cy3structure} one has
\begin{equation}\label{eq:yukawa}
\langle e_{2h+1},
\nabla\left(\frac{\partial^3 }{\partial \tau_i\partial \tau_j\partial
\tau_k}\right)e_{2h+1}\rangle\;=\;
-\frac{1}{2}\frac{\partial^3 Z}{\partial \tau_i\partial \tau_j\partial \tau_k}.
\end{equation}
here
$\nabla\left(\frac{\partial^3 }{\partial \tau_i\partial \tau_j\partial
\tau_k}\right)$ is shorthand for $ \left(\frac{\partial}{\partial
\tau_i}\otimes 1\right)
\circ\nabla\circ\left(\frac{\partial}{\partial \tau_j}\otimes 1\right)
\circ\nabla\circ\left(\frac{\partial}{\partial \tau_k}\otimes
1\right)\circ\nabla$.
\end{corollary}
\begin{proof}
\begin{eqnarray*}
&&\hspace{-3em}\langle e_{2h+1},\nabla\left(\frac{\partial^3 }{\partial
\tau_i\partial \tau_j\partial \tau_k}\right)
e_{2h+1}\rangle\;=\;
\frac{\partial^3 \tau_{03}}{\partial \tau_i\partial \tau_j\partial
\tau_k}-\tau_{23}^\star\frac{\partial^3 \tau_{13}}{\partial \tau_i\partial
\tau_j\partial \tau_k}\,=
\\&&\hspace{2em}=\;
\frac{\partial^3 Z}{\partial \tau_i\partial \tau_j\partial \tau_k}
-\frac{1}{2}\frac{\partial^2}{\partial \tau_j\partial \tau_k}
\frac{\partial Z}{\partial \tau_i}
-\frac{1}{2}\frac{\partial^2}{\partial \tau_i\partial \tau_k}
\frac{\partial Z}{\partial \tau_j}
-\frac{1}{2}\frac{\partial^2}{\partial \tau_i\partial \tau_j}
\frac{\partial Z}{\partial \tau_k}.
\end{eqnarray*}
\end{proof}

\

All hypotheses and notations of Theorem \ref{cy3structure} remain in force.
The divisibility hypothesis (\ref{eq:divisible}) and Formula
(\ref{eq:frobmat2}) imply that for every lift of Frobenius $\psi$ the matrix
$\cT^{-1}\cdot P\cdot \psi(\cT)\cdot P^{-1}$ has entries in $\oA$.
This means in particular
$$
p^{-1}\psi(\tau_i)-\tau_i\in\oA\qquad\textrm{for}\quad i=1,\ldots,h.
$$
{}From Theorem \ref{thmstructure} we also know that
$\tau_i(0)\in p\oW$ for $i=1,\ldots,h$.
By a lemma of Dwork this implies (see \cite{D1} cor.1.4.5) that, if we set
$$
q_i=exp(\tau_i)\qquad\textrm{for}\quad i=1,\ldots,h,
$$
then $q_i$ is well defined and
\begin{equation}\label{eq:canorig}
q_i\in\oA\qquad\textrm{and}\qquad q_i(0)-1\in p\oW.
\end{equation}
\cite{D1} cor.1.4.7 now yields that $p,q_1-1,\ldots,q_h-1$ is a regular system
of parameters for the ring $\oA$ i.e.
$$
\oA=\oW[[t_1,\ldots,t_h]]\:=\:\oW[[q_1-1,\ldots,q_h-1]].
$$
The elements $q_1,\ldots,q_h$ are called the \emph{canonical coordinates}.
There is an associated \emph{canonical lift of Frobenius}
$$
\vf:\oA\rightarrow\oA\,,\qquad \vf|_\oW=\sigma\,,\qquad \vf(q_i)=q_i^p
\quad\textrm{for}\;i=1,\ldots,h.
$$

\

For every lift of Frobenius $\psi$ the matrix of the map $F(\psi)\psi^*$ with
respect to the basis $e_0,\ldots,e_{2h+1}$ is given by Formula
(\ref{eq:frobmat2}).
A straightforward computation makes the entries of this matrix explicit.
One then sees that exactly for the canonical lift of Frobenius $\vf$ this
matrix takes the following elegant form:
\\
\emph{
For the canonical lift of Frobenius $\vf$
the matrix of $F(\vf)\vf^*$ with respect to the basis $e_0,\ldots,e_{2h+1}$ is}
\begin{equation}\label{eq:matcanfrob}
\left(\begin{array}{cccc}
1&0&-p^{-2}\vf(\tau_{13}^\star)+\tau_{13}^\star&p^{-3}\vf(Z)-Z\\
0&I_h&p^{-1}\vf(\tau_{12})-\tau_{12}&p^{-2}\vf(\tau_{13})-\tau_{13}\\
0&0&I_h&0\\
0&0&0&1
\end{array}\right)
\left(\begin{array}{cccc}
1&0&0&0\\ 0&pI_h&0&0\\0&0&p^2I_h&0\\ 0&0&0&p^3
\end{array}\right)
\end{equation}
This formula
generalizes \cite{D1} (1.4.7.1) to the CY3 case.

We summarize the preceding discussion in
\begin{theorem}\label{integralities}
All hypotheses and notations of Theorem \ref{cy3structure} remain in force.
Define
\begin{equation}\label{eq:qcoord}
q_i=exp(\tau_i)\qquad\textrm{for}\quad i=1,\ldots,h.
\end{equation}
Then
\begin{equation}\label{eq:qprops}
q_i\in\oA\qquad\textrm{and}\qquad q_i(0)-1\in p\oW.
\end{equation}
and $p,q_1-1,\ldots,q_h-1$ is a regular system of parameters for the ring $\oA$
i.e.
\begin{equation}\label{eq:ttoq}
\oA=\oW[[t_1,\ldots,t_h]]\:=\:\oW[[q_1-1,\ldots,q_h-1]].
\end{equation}
The elements $q_1,\ldots,q_h$ are called the \emph{canonical coordinates}.
There is an associated \emph{canonical lift of Frobenius}
\begin{equation}\label{eq:canonicalfrob}
\vf:\oA\rightarrow\oA\,,\qquad \vf|_\oW=\sigma\,,\qquad \vf(q_i)=q_i^p
\quad\textrm{for}\;i=1,\ldots,h.
\end{equation}
Moreover
\begin{equation}\label{eq:yukinnerproduct}
\langle e_{2h+1},F(\vf)\vf^*(e_{2h+1})\rangle=\vf(Z)-p^3Z
\end{equation}
and
\begin{equation}\label{eq:yukinteger}
p^{-3}\vf(Z)-Z\in\oA=\oW[[t_1,\ldots,t_h]].
\end{equation}
\end{theorem}
\begin{proof}
(\ref{eq:yukinnerproduct}) and (\ref{eq:yukinteger}) follow directly from
(\ref{eq:matcanfrob}) and the divisibility hypothesis (\ref{eq:divisible}).
\end{proof}

\section{Ordinary CY3-crystals over a perfect field}\label{sectionCY3p}
Having a good hold on the structures when the base field is algebraically
closed we now turn to a base field which is just a \emph{perfect field} $k$.
So:\\
\emph{Let $p>3$ and let $H$ be a divisible ordinary Hodge $F$-crystal of level
$3$ over $A_0$ equipped with a non-degenerate alternating bilinear form
$\langle\:,\:\rangle$
such that (\ref{eq:Hodgenumbers}), (\ref{eq:formconditions}),
(\ref{eq:Riemann}) hold and such that (\ref{eq:conngrhodge}) is an
isomorphism}, everything with $H$ in place of $\oH$.

Let $\oH$ be the Hodge $F$-crystal over $\oA_0$ obtained from $H$ by base
change from $k$ to $\ok$. Then all results of Section \ref{sectionCY3} hold for
$\oH$.
We use the notations of  Section \ref{sectionCY3} for referring to matters of
$\oH$.

Fix a non zero element $\omega_0$ in $\Fil^3 H=H^3$. Since the $\oA$-module
$\Fil^3 \oH=\oH^3$ has rank $1$ with basis $e_{2h+1}$, there is an invertible
element $f\in\oA$ such that
\begin{equation}\label{eq:omegaf}
\omega_0=f e_{2h+1}.
\end{equation}

\begin{lemma}\label{lemmaf}
The element $f$ in (\ref{eq:omegaf}) has the form
\begin{equation}\label{eq:f}
f=a\tf^{-1}\qquad\textrm{with}\quad a\in\oW,\;\tf\in A,\;\tf(0)=1.
\end{equation}
\end{lemma}
\begin{proof}
Let $\hat{\omega}\in V_3$ and $\hat{e}_{2h+1}\in\oV_3$ denote the respective
images
of $\omega_0$ and $e_{2h+1}$ under the isomorphisms
$
\Fil^3H\stackrel{\simeq}{\rightarrow}U_3/U_2 = V_3\;,\quad
\Fil^3\oH\stackrel{\simeq}{\rightarrow}\oU_3/\oU_2 = \oV_3.
$
Take the lift of Frobenius $\phi:A\rightarrow A$ such that $\phi (t_i)=t_i^p$
for $i=1,\ldots,h$. Then in the unit $F$-crystal $V_3$
$$
F(\phi)\phi^*(\hat{\omega})= c g \hat{\omega}
$$
with $c\in W$, $g\in A$, $g(0)=1$.
The product
$
\tf=\prod_{i=0}^\infty \phi^i(g)
$
converges in $A$, because $\phi^i(g)\equiv 1\bmod
(t_1^{p^i},\ldots,t_h^{p^i})$.
Moreover there is an element $b\in\oW$ such that
$c=\sigma(b)b^{-1}$.
Thus
$$
F(\phi)\phi^*(b^{-1}\tf\hat{\omega})=
\sigma(b^{-1})\phi(\tf)cg\hat{\omega}=
b^{-1}\tf\hat{\omega}.
$$
On the other hand, (\ref{eq:frobmat2}) implies
$F(\phi)\phi^*e_{2h+1}=p^3e_{2h+1}$
and hence in the (untwisted) unit $F$-crystal $\oV_3$ we have
$$
F(\phi)\phi^*\hat{e}_{2h+1}=\hat{e}_{2h+1}.
$$
Since
$
\hat{\omega}=f\hat{e}_{2h+1}
$
we see that
$\phi(b^{-1}\tf f)=b^{-1}\tf f$ and hence
$b^{-1}\tf f\in\ZZ_p$. Thus $f=a \tf^{-1}$ with $a\in \oW$.
\end{proof}

\

We set for $j=1,\ldots,h$
\begin{equation}\label{eq:omj1}
\omega_j:=
fD_j\left(\frac{1 }{f}\omega_0\right)
\end{equation}
where
$D_j:=(\frac{d}{dt_j}\otimes 1)\circ\nabla$.
Lemma \ref{lemmaf} shows
$\omega_j=
\tf^{-1} D_j\left(\tf\omega_0\right)$
and hence $\omega_j\in H$.
\\
On the other hand, $\nabla e_{2h+1}=\sum_{i=1}^h d\tau_i\otimes e_{h+i}$
implies
\begin{equation}\label{eq:omj2}
\omega_j=fD_j(e_{2h+1})=f\frac{\partial \tau_1}{\partial t_j}e_{h+1}+\ldots+
f\frac{\partial \tau_h}{\partial t_j}e_{2h}.
\end{equation}
{}From (\ref{eq:kodaira-spencer}) we know that
$\left(\frac{\partial\tau_i}{\partial t_j}\right)_{1\leq i,j\leq h}$ is
precisely the matrix of the isomorphism
$T_{\oA/\oW}\rightarrow \mathrm{Hom}(\oH^3,\oH^2)$.
Since $\{e_{h+1},\ldots,e_{2h}\}$ is an $\oA$-basis of $\oH^2$, we now see that
$\{\omega_1,\ldots,\omega_h\}$ is an $A$-basis for $H^2$.

Next we set for $i=1,\ldots,h$
\begin{equation}\label{eq:omj3}
\check{\omega}_i=f^{-1}\frac{\partial t_i}{\partial \tau_1}e_1+\ldots+
f^{-1}\frac{\partial t_i}{\partial \tau_h}e_h.
\end{equation}
Then
\begin{equation}\label{eq:omj4}
\langle \check{\omega}_i,\omega_j\rangle\:=\:\sum_{m=1}^h
\frac{\partial t_i}{\partial \tau_m}
\frac{\partial \tau_m}{\partial t_j}\:=\:
\left\{\begin{array}{lll}
0&\textrm{ for }&i\neq j\\
1&\textrm{ for }&i= j
\end{array}\right..
\end{equation}
The form $\langle\,,\,\rangle$ induces a duality between $H^1$ and $H^2$. This
fact together with (\ref{eq:omj4}) shows that
$\check{\omega}_1,\ldots,\check{\omega}_h$ which a priori are in
$\oH^1$, lie in fact in $H^1$ and constitute the basis of $H^1$ dual to the
basis $\omega_1,\ldots,\omega_h$ of $H^2$.

Finally we set
\begin{equation}\label{eq:omj5}
\check{\omega}_0=f^{-1}e_0.
\end{equation}
Then $\check{\omega}_0\in \oH^0$ and $\langle
\check{\omega}_0,\omega_0\rangle=-1$.
So $\check{\omega}_0$ lies in fact in $H^0$ and gives a basis dual to the basis
$\omega_0$ of $H^3$.

We summarize the preceding discussion in the following proposition.

\begin{proposition}\label{propomega}
Starting from $\omega_0\in\Fil^3H=H^3$ define $f\in\oA$ by (\ref{eq:omegaf}),
define $\omega_1,\ldots,\omega_h\in H^2$  by (\ref{eq:omj1}) and
define $\check{\omega}_0\in H^0$ and
$\check{\omega}_1,\ldots,\check{\omega}_h\in H^1$ by
(\ref{eq:omj5}) resp. (\ref{eq:omj3}).
Write
$$
\check{\omega}_{2h+1}=\omega_0\,,\qquad
\check{\omega}_{h+i}=\omega_i\quad\textrm{for}\;i=1,\ldots,h.
$$
Then $\{\check{\omega}_0,\ldots,\check{\omega}_{2h+1}\}$ is a basis
for $H$ such that
$$
\check{\omega}_0\in H^0\,,\quad \check{\omega}_1,\ldots,\check{\omega}_h\in
H^1\,,\quad \check{\omega}_{h+1},\ldots,\check{\omega}_{2h}\in H^2
\,,\quad \check{\omega}_{2h+1}\in H^3.
$$
The matrix $\cS$ in the basis transformation
$$
(\check{\omega}_0,\ldots,\check{\omega}_{2h+1})=(e_0,\ldots,e_{2h+1})\cS
$$
(i.e. the $j$-th column gives the coordinates of $\check{\omega}_j$ with
respect to $e_0,\ldots,e_{2h+1}$) is
\begin{equation}\label{eq:matomep}
\cS=
\left(
\begin{array}{cccc}
f^{-1}&0&0&0\\
0&f^{-1}
\left(\frac{\partial t_j}{\partial\tau_i}\right)
&0&0\\
0&0&f
\left(\frac{\partial\tau_i}{\partial t_j}\right)&0\\
0&0&0&f
\end{array}
\right)
\end{equation}
\qed
\end{proposition}

Recall that in the proof of Theorem \ref{cy3structure} we constructed a basis
$\{\ep_0,\ldots,\ep_{2h+1}\}$
from the basis transformation
$$
(e_0,\ldots,e_{2h+1})=(\ep_0,\ldots,\ep_{2h+1})\cT.
$$
The relation with the basis $\{\check{\omega}_0,\ldots,\check{\omega}_{2h+1}\}$
is therefore
$$
(\check{\omega}_0,\ldots,\check{\omega}_{2h+1})=
(\ep_0,\ldots,\ep_{2h+1})\cT\cS.
$$
The virtues of $\ep_0,\ldots,\ep_{2h+1}$ are
$$
\nabla\ep_i=0 \qquad\textrm{for}\quad i=0,\ldots,2h+1
$$
and for every lift of Frobenius $\phi:A\rightarrow A$
\begin{eqnarray*}
&&F(\phi)\phi^*\ep_0=\ep_0\,,\quad F(\phi)\phi^*\ep_{2h+1}=p^3\ep_{2h+1}\\
&&F(\phi)\phi^*\ep_i=p\ep_i\,,\quad F(\phi)\phi^*\ep_{h+i}=p^2\ep_{h+i}
\qquad\textrm{for}\quad i=1,\ldots,h.
\end{eqnarray*}
Thus one finds the following analogue of formulas (\ref{eq:conmat2}) and
(\ref{eq:frobmat2})

\begin{corollary}
Let the hypotheses and notations be as in Proposition \ref{propomega}.
Then the matrix of the connection $\nabla$ with respect to the basis
$\{\check{\omega}_0,\ldots,\check{\omega}_{2h+1}\}$ is
\begin{equation}\label{eq:conmatomega}
(\cT\cS)^{-1}\cdot d(\cT\cS)
\end{equation}
and for every lift of Frobenius $\phi:A\rightarrow A$ the matrix of the map
$F(\phi)\phi^*$ with respect to the basis
$\{\check{\omega}_0,\ldots,\check{\omega}_{2h+1}\}$ is
\begin{equation}\label{eq:frobmatomega}
(\cT\cS)^{-1}\cdot P\cdot\phi(\cT\cS)
\end{equation}
with matrices $\cT$ and $P$ as in Theorem \ref{cy3structure} and $\cS$ as in
Proposition \ref{propomega}.
\qed
\end{corollary}

\

\begin{remark}\label{Picard-Fuchs}
Up to now we have discussed all structures in terms of a connection, to be
thought of as the Gauss-Manin connection of a family of Calabi-Yau threefolds.
In the literature descriptions of canonical coordinates and Yukawa coupling for
complex Calabi-Yau threefolds near the large complex structure limit usually
start from the Picard-Fuchs equations and their solutions.
We want to point out how such a description can also be seen in our discusssion
of ordinary CY3 crystals.
The hypotheses and notations are as before in this section.

We start with a non zero element $\omega_0$ in $\Fil^3 H$. In the geometric
situation $\Fil^3 H$ would be $H^{3,0}$ and thus $\omega_0$ is the analogue of
a nowhere vanishing global $3$-form.
The connection induces an action of differential operators (in
$\frac{\partial}{\partial t_1},\ldots,\frac{\partial}{\partial t_h}$)
on $H$. The differential operators which annihilate $\omega_0$ are called
\emph{the Picard-Fuchs operators}.
We expand $\omega_0$ with respect to the basis $\{\ep_0,\ldots,\ep_{2h+1}\}$:
$$
\omega_0=f_0\ep_0+\ldots+f_{2h+1}\ep_{2h+1}.
$$
Since $\nabla\ep_i=0$ for all $i$, all $f_i$ are annihilated by
the Picard-Fuchs operators, i.e. are solutions of the Picard-Fuchs equations.
This set of solutions is split into subsets
$\{f_0\}\,,\;\{f_1,\ldots,f_h\}\,,\;\{f_{h+1},\ldots,f_{2h}\}\,,\;
\{f_{2h+1}\}$ according to the action of Frobenius on the basis vectors
$\ep_0,\ldots,\ep_{2h+1}$. In the complex setting near the large complex
structure point such a splitting of a basis of the solution space of the
Picard-Fuchs equations is made according to the monodromy action, i.e.
according to the degree of the logarithmic terms in the solution.

Since $\check{\omega}_{2h+1}=\omega_0$ and
$
(\check{\omega}_0,\ldots,\check{\omega}_{2h+1})=
(\ep_0,\ldots,\ep_{2h+1})\cT\cS
$
the column vector $(f_0,\ldots,f_{2h+1})^\star$ is the last column of the
matrix $\cT\cS$. Thus $f=f_{2h+1}$ and $\left(\frac{f_0}{f},\ldots,
\frac{f_{2h+1}}{f}\right)$ is the last column of the matrix $\cT$.
This is analogous to the step of dividing the solutions of the Picard-Fuchs
equations by the unique solution which is holomorphic and has value $1$
at the large complex structure point.

{}From here on the algorithms for computing the canonical coordinates and the
prepotential of the Yukawa coupling are identical for the ordinary case and for
the large complex structure case.
\end{remark}
%%%%%%%%%%%%%%%%%

%%
\


\begin{thebibliography}{99}
%
%
\bibitem{BK}
Bloch, S., K. Kato, \textit{$p$-Adic \'etale cohomology}, Publ. Math. IHES
\textbf{63}, 1986, 107--152.
%
\bibitem{BG}
Bryant, R., and P. Griffiths, \textit{Some observations on the infinitesimal
period relations for regular threefolds with trivial canonical bundle},
Arithmetic and Geometry vol. II, M. Artin and J. Tate (eds.), Progress in Math.
vol. 36, Birkh\"auser, Boston, 1983, pp. 77--102.
%
\bibitem{C}
Candelas, P., X. de la Ossa, P. Green, L. Parkes,
\textit{A pair of Calabi-Yau manifolds as an exactly soluble superconformal
theory}, Essays on Mirror Manifolds, S-T. Yau (ed.),International Press, Hong
Kong, 1992, pp. 31--95.
%
\bibitem{D1}
Deligne, P., \textit{Cristaux ordinaires et coordonn\'ees canoniques},
Surfaces Alg\'ebriques, J. Giraud, L. Illusie and M. Raynaud (eds.), Lecture
Notes in Math. 868, Springer-Verlag, Berlin, 1981, pp. 80--127.
%
\bibitem{D2}
Deligne, P., \textit{Local behavior of Hodge structures at infinity}, Mirror
Symmetry II, B. Greene and S-T. Yau (eds.), AMS/IP Studies in Advanced Math.
vol. 1, Amer. Math. Soc. and International Press, 1997, pp. 683--699.
%
\bibitem{F}
Friedman, R., \textit{On threefolds with trivial canonical bundle},
Complex Geometry and Lie Theory, J. Carlson, H. Clemens and D. Morrison (eds.),
Proc. Symp. Pure Math. vol.53, Amer. Math. Soc., Providence, RI, 1991,
pp.103--134.
%
\bibitem{K1}
Katz, N., \textit{Travaux de Dwork}, S\'eminaire Bourbaki exp.409, Lecture
Notes in Math. 317, Springer-Verlag, Berlin, 1973, pp.167--200.
%
\bibitem{M}
Morrison, D., \textit{Mirror symmetry and rational curves on quintic
threefolds: a guide for mathematicians}, J. of the Amer. Math. Soc. \textbf{6}
(1993), 223--247; alg-geom/9202004
%
\bibitem{LY}
Lian, B., and S-T. Yau, \textit{Mirror maps, modular relations and
hypergeometric series I},  hep-th/9507151
%
\bibitem{S}
Stienstra, J., \textit{The ordinary limit for varieties over
$\ZZ[x_1,\ldots,x_r]$}, these proceedings.
\end{thebibliography}
\end{document}